\documentclass{amsart}

\usepackage{graphicx}
\usepackage{amsfonts}
\usepackage{amsmath}
\usepackage{xcolor}
\usepackage{multicol}
\usepackage{caption}
\usepackage{subcaption}

\DeclareMathOperator*{\argmax}{argmax}  
\DeclareMathOperator*{\argmin}{argmin}

\newtheorem{theorem}{Theorem}[section]
\newtheorem{lemma}[theorem]{Lemma}

\theoremstyle{definition}
\newtheorem{definition}[theorem]{Definition}

\newtheorem{corollary}[theorem]{Corollary}

\numberwithin{equation}{section}

\begin{document}

\title[The SCFE Problem in the Circumference]{The Sitting Closer to Friends than Enemies Problem in the Circumference}

\author{Julio Aracena}\address{CI$^2$MA \and Departamento de Ingenier\'ia Matem\'atica, Facultad de Ciencias F\'isicas y Matem\'aticas, Universidad de Concepci\'on, Chile.}\email{jaracena@ing-mat.udec.cl} \thanks{Julio Aracena was partially supported by ANID-Chile through the project {\sc Centro de Modelamiento Matem\'atico} (AFB170001) of the PIA Program:``Concurso Apoyo a Centros Cient\'ificos y Tecnol\'ogicos de Excelencia con Financiamiento Basal'', and by ANID-Chile through Fondecyt project 1151265.}

\author{Felipe Benitez}\address{Department of Mathematics and Statistics, Georgia State University, Atlanta, USA.}\email{fbenitez1@student.gsu.edu}\thanks{Felipe Benitez was partially supported by CONICYT-Chile through Fondecyt project 1151265.}

\author{Christopher Thraves Caro}\address{Departamento de Ingenier\'ia Matem\'atica, Facultad de Ciencias F\'isicas y Matem\'aticas, Universidad de Concepci\'on, Chile.}\email{cthraves@ing-mat.udec.cl}

%    General info
\subjclass[2000]{Primary 05C22, 05C62, 68R10. Secondary 05C85}  

\date{}

\keywords{Signed Graphs, The SCFE problem, Valid Distance Drawing, Proper Circular-Arc Graphs, Circular-Arc Graphs} %TODO mandatory; please add comma-separated list of keywords

% REQUIRED
\begin{abstract}
Consider a graph $G=(V, E)$ and a sign assignment to the edges of $G$. The \emph{Sitting Closer to Friends than Enemies} (SCFE) problem is to find an injection of $V$ into a metric space so that, for every pair of incident edges with a different sign, the vertices positively connected are closer (in the metric of the space) than the vertices negatively connected. In this document, we present recent results regarding the SCFE problem when the metric space in consideration is the circumference. In particular, we prove that given a complete signed graph, it has an injection that satisfies the SCFE problem in the circumference if and only if its positive subgraph is a proper circular-arc graph. 
\end{abstract}

\maketitle
\section{Introduction}
\label{intro}
Consider a group of people that may have positive or negative interactions between them. For instance, they may be friends or enemies, or they may not know each other. Now, sit them all at a big circular table so that their friends first surround each person and then, farther than all their friends, the person can see their enemies. Finding such a placement in a circular table is known as the \emph{Sitting Closer to Friends than Enemies} (SCFE) problem in the circumference. 

A way to represent the group of people is using a \emph{signed graph}, i. e., a graph with an assignment of signs, positive or negative, on its edges. Therefore, each vertex of the graph represents one person. We connect two friends with a positive edge and two enemies with a negative edge. Two unknown people are not adjacent in the graph. In this document, we show that when the signed graph is complete (all its edges are present), we can solve the SCFE problem in the circumference in polynomial time. Furthermore, we characterize the set of complete signed graphs for which the SCFE problem in the circumference has a solution. 

\section{Definitions and Notation}
In this manuscript, we consider signed graphs that are finite, 
undirected, connected, loopless and without parallel edges.
We use $S=(V,E^+\cup E^-)$ to denote a signed graph with vertex set $V(S)$, set of positive edges $E^+(S)$, and set of negative edges
$E^-(S)$.
When the signed graph under consideration is clear, we use $V$, 
$E^+$ and $E^-$. It is worth noting that in every signed graph $E^+\cap 
E^-=\emptyset$.  
The number of vertices of a signed graph $S$ is denoted by $n :=|V(G)|$,
the number of positive edges is denoted by $m^+:=|E^+(G)|$, and the total
number of edges of a graph
(positive plus negative edges) is denoted by $m:=|E^+(G)|+|E^-(G)|$. 

A positive or negative edge $\{i,j\}$ is denoted $ij$. If $ij\in E^+$ 
(resp., $ij \in E^-$), we say that $i$ is a \emph{positive neighbor} 
(resp., \emph{negative neighbor}) of $j$ and vice versa. 
The set of positive (resp., negative) neighbors of $i$ is denoted 
$N^+(i)$  (resp., $N^-(i)$). Additionally, the \emph{closed positive 
neighborhood} of $i$ is defined as $N^+[i]:= N^+(i) \cup \{i\}.$

A signed graph $H$ is a \emph{signed subgraph} of a signed graph $S$ 
if and only if: 
\[V(H) \subseteq V(S) \quad \land \quad E^+(H) \subseteq E^+(S) \quad 
\land \quad E^-(H) \subseteq E^-(S).\]
The \emph{positive subgraph} of a signed graph $S$ is the signed graph $S^+ =(V(S), E^+(S) \cup \emptyset)$ that contains all 
the vertices, all the positive edges, and none of the negative 
edges of $S$. Even though, the positive subgraph of a signed graph is
a signed graph, it can also be seen as a graph (with no signs on its 
edges) since all its edges have the same sign.
We will use $G$ to denote graphs with no signs.  
A signed graph is \emph{complete} if, for every 
pair of distinct vertices $i$ and $j$, $ij \in E^+$ or $ij \in E^-$.

The \emph{circumference}, w.l.o.g.\footnote{It is worth noting that this
definition can be modified 
changing the radius and/or the center of the circumference and all the 
results presented in this document remain valid.}, 
is defined as the set $\mathcal{C} = \{(x,y):\sqrt{x^2+y^2}=1\} \subseteq
\mathbb{R}^2$.  Every point $p$ in $\mathcal{C}$ is fully determined by 
the angle, in $[0,2\pi[$, formed by the points $(1,0)$, the origin 
$\vartheta=(0,0)$ and the point $p$ (moving clockwise). Hence, as
an abuse of notation, we use that angle in $[0,2\pi[$ to denote $p \in 
\mathcal{C}$. For example, the point $(1,0)$ is the point $0$, the 
point $(0,-1)$ is the point $\pi/2$, and the point $(0,1)$ is the point $3\pi / 2$. The \emph{distance} between two 
points $p$ and $q$ in $\mathcal{C}$ is the size of the smallest angle 
formed by $p$, $\vartheta$, and $q$. Hence, the distance 
between $p$ and $q$ is defined as: 
\[
d(p,q) := \min \{{(p-q)}\bmod {2\pi},{(q-p)}\bmod {2\pi} \}.
\]
The pair $(\mathcal{C},d)$ is the metric space that we consider in this 
work. As an abuse of notation, we will refer to $(\mathcal{C},d)$ as $\mathcal{C}$.

A \emph{drawing} of a signed graph $S$ in $\mathcal{C}$ is an injection 
$D:V\rightarrow \mathcal{C}$ of the vertex set $V$ into $\mathcal{C}$. Given a drawing $D$ of $S$ in $\mathcal{C}$, and $i\in V$, the function $d_{i}^{D}:V\rightarrow [0,2\pi[$ defined by $d_{i}^{D}(j)=d(D(i),D(j))$, $\forall j\in V$, represents the distance from $D(i)$ to $D(j)$.
A drawing $D$ of $S$ in $\mathcal{C}$ is said to be \emph{valid distance} if, for all
$i \in V$, for all $j \in N^+(i)$, and for all  $k \in N^-(i)$,
\begin{equation}\label{cond:valid} 
d_{i}^{D}(j) < d_{i}^{D}(k).
\end{equation}

\begin{figure}[t]
\begin{minipage}[t]{.33\linewidth}
 \centering\includegraphics[scale=1]{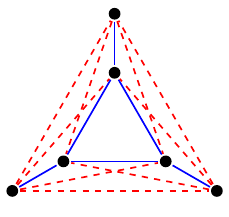}
 \subcaption{Signed graph without a valid distance drawing}\label{fig:2:1}
 \end{minipage}
 \begin{minipage}[t]{.32\linewidth}
 \centering\includegraphics[scale=1]{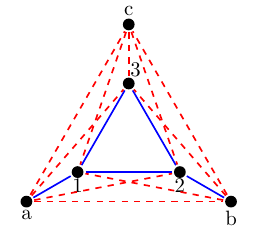}
 \subcaption{Signed graph with a valid distance drawing}\label{fig:2:2}
 \end{minipage}
 \begin{minipage}[t]{.32\linewidth}
 \centering\includegraphics[scale=1]{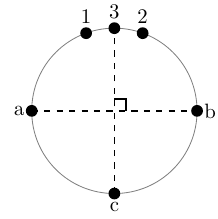}
 \subcaption{Valid distance drawing of the signed graph in Subfigure \ref{fig:2:2}}\label{fig:2:3}
 \end{minipage}
 \caption{Subfigure \ref{fig:2:1} and Subfigure \ref{fig:2:2} show two different complete signed 
graphs where dashed lines represent negative edges and continuous lines 
represent positive edges. The positive subgraphs of these signed graphs 
are the subgraphs composed of the continuous edges only.  The positive 
subgraph of the signed graph in Subfigure \ref{fig:2:1} 
is known as the \emph{net}. The net is not a proper ciruclar-arc graph. As we will show in Theorem \ref{thm:characcomplete},
the complete signed graph in Subfigure \ref{fig:2:1} does not have a valid distance drawing since the net is not a proper circular-arc graph. The complete signed graph in Subfigure \ref{fig:2:2} has a valid distance drawing and Subfigure \ref{fig:2:3} shows one.\label{fig}}
 \end{figure}

In the case that there exists a valid distance drawing of a given signed graph $S$
in $\mathcal{C}$, we say that $S$
\emph{has a valid distance drawing} in $\mathcal{C}$. Otherwise, we say 
that $S$ is a signed graph without valid distance drawing in $\mathcal{C}$. 
The definition of a valid distance drawing 
captures the requirement that every 
vertex is placed closer to its friends than to its enemies.
Hence, the Sitting Closer to Friends than Enemies problem in the 
circumference is: 
\begin{definition}{\emph{(SCFE problem in the circumference)}}
Given a signed graph $S$ decide whether $S$ has a valid distance 
drawing in 
$\mathcal{C}$.
\end{definition}
If, for a given signed graph $S$, the SCFE problem has a positive answer,
then, we are interested in finding a valid distance drawing for $S$ in 
$\mathcal{C}$. 
Figure \ref{fig} shows a complete signed graph without a valid distance drawing in 
$\mathcal{C}$, a
complete signed graph with a valid distance drawing in 
$\mathcal{C}$, and a valid distance drawing of that graph in 
$\mathcal{C}$.

%\begin{figure}[t]
%\begin{center}
%\resizebox{.8\textwidth}{!}{\includegraphics*{Dibujo3grafos.pdf}}
%\caption{Figure (a) and Figure (b) show two different complete signed 
%graphs where dashed lines represent negative edges and continuous lines 
%represent positive edges. The positive subgraphs of these signed graphs 
%are the subgraphs composed of the continuous edges only.  The positive 
%subgraph of the signed graph in Figure (a) is known as the \emph{net}. 
%The complete signed graph in Figure (a) does not have a valid distance %drawing since the net is not a proper circular-arc graph. The complete %signed graph in Figure (b) has a valid distance drawing and Figure (c) %shows one.}
%\label{fig}
%\end{center}
%\end{figure}

Given $p \in \mathcal{C}$, we 
define the \emph{right half} and \emph{left half} of $p$ in $\mathcal{C}$
as $M_r(p):=\{(p+t)\bmod2\pi:0\leq t\leq \pi\}$  and 
$M_l(p):=\{(p-t)\bmod2\pi:0\leq t\leq \pi\}$, respectively. 
Therefore, given two points $p$, and $q$ in  $\mathcal{C}$, 
\[p \in M_r(q) \iff 
q\in M_l(p).\]  
It is worth mentioning that, for every point 
$p \in \mathcal{C}$, $p\in M_l(p)$ and $p \in M_r(p)$.

Finally, a drawing of a signed graph in $\mathcal{C}$ induces a cyclic 
order of its vertices. Given a signed graph $S$ and a drawing $D$ of 
$S$ in $\mathcal{C}$, 
we say that $i$ is smaller than $j$ according to $D$ if, starting from 
the point $0$ and traveling $\mathcal{C}$ in clockwise direction, we find first $i$ and then $j$. 
In such case, we denote $i<_Dj$ or simply $i<j$ if the drawing $D$ is 
clear by the context. Now, given a signed graph $S$ and a drawing $D$ of
$S$ in $\mathcal{C}$, we 
relabel the vertices naming $0$ the first vertex in the ordering induced
by $D$, $1$ the second vertex, and so on until $n-1$, the last vertex. 
Hence,  $0 <_D 1<_D 2<_D\cdots <_D n-1$. 
It is worth noting that it also holds $D(0) < D(1) < \cdots <D(n-1)$, when,
again, $D(i)$, as an abuse of notation, is the angle in $[0,2\pi[$ formed
by the point $0$,  the origin $\vartheta$, and the 
point in which $D$ injects vertex $i$ (clockwise). Even though, the order as defined
above is not cyclic, we provide it with a cyclic structure due to the 
circular characteristic of the space. Hence, the vertex set
$V=\{0,1,2,\ldots,n-1\}$ is considered cyclically ordered, i.e., 
$0<1<2<\ldots < n-1 <0$. 

\section{Related Work and Our Contributions}
The SCFE problem was first introduced by Kermarrec et al. in \cite{kermarrec2011can}. In their work, Kermarrec and Thraves studied the SCFE problem in $\mathbb{R}$. 
They presented families of signed graphs without a 
valid distance drawing in 
$\mathbb{R}$. They also 
gave a characterization of the set of signed graphs 
with a valid distance drawing in $\mathbb{R}$.
Such characterization was then used by the authors to construct 
a polynomial time decision algorithm for complete signed graphs with a valid distance drawing in $\mathbb{R}$. 

Afterwards, Cygan et al. in \cite{cygan2015sitting} proved that the SCFE
problem in $\mathbb{R}$ is NP-Complete. In addition, they gave another
characterization of the set of complete signed graphs with a 
valid distance drawing in $\mathbb{R}$. This characterization says that a
complete signed graph $S$ 
has a valid distance drawing in $\mathbb{R}$ if and only if $S^+$ is a proper
interval graph. 

Based on the previous NP-Completeness result, Garcia Pardo et al. in 
\cite{pardo2015embedding} studied an optimization version of the SCFE 
problem in $\mathbb{R}$ where the goal is to find a drawing that 
minimizes the
number of violations of condition \eqref{cond:valid} in the definition of a valid distance 
drawing. They proved that when the signed graph $S$ is complete,
local minimums for their optimization problem 
coincide with local minimums of the well known \emph{Quadratic 
Assignment} problem
applied to $S^+$. Moreover, they studied experimentally two heuristics, 
showing that a greedy heuristic has a good performance at the moment of 
recognition of graphs that have an optimal solution with zero errors 
(problem that is NP-Complete). Pardo et al. in \cite{pardo2020basic}
improved the experimental results on this problem presenting a \emph{basic 
variable neighborhood search} algorithm.

Spaen et al. in \cite{spaen2020dimension} studied the SCFE problem from a 
different perspective. They studied the problem of finding
$L(n)$, defined as the smallest dimension $k$ such that any signed graph on $n$
vertices has a valid distance drawing in $\mathbb{R}^k$, with respect
to Euclidean distance. They showed that $\log_5 (n-3) \leq L(n) \leq n - 
2$.

Finally, Becerra et al. in \cite{becerra2019sitting} studied the SCFE problem 
in trees. They proved that a complete signed 
graph $S$ has a valid distance drawing in a tree if and only if $S^+$ is 
strongly chordal. 

\paragraph*{Our Contributions}
In this document, we present our study on the SCFE problem in the circumference. The main result of this work is presented in Theorem \ref{thm:characcomplete} and says that a complete signed graph $S$ has a valid distance drawing in $\mathcal{C}$ if and only if $S^+$ is a proper circular-arc graph. The rest of the document has the following structure. In Section \ref{sec:instrumental}, we present two instrumental lemmas that we use in Section \ref{sec:characterization} to prove our main result. Finally, in Section \ref{sec:concludingremarks}, we present some concluding remarks. 

\section{Instrumental Results}\label{sec:instrumental}
In this section, we present two  results that we use in the rest of the document. 
\subsection{Almost Valid Distance Drawings}\label{sec:instlemmas}
 With instrumental purposes, we introduce a different type of drawing that we call \emph{almost valid distance}. 
 Let $S=(V,E^+ \cup E^-)$ be a signed graph and $D$ be a drawing of $S$ in $\mathcal{C}$. 
 We say that $D$ is \emph{almost valid distance} if and only if  there exists  $\delta>0$ such that for all
 $i \in V$, for all $j \in N^+(i)$ and for all $k \in N^-(i)$,
\begin{equation}\label{cond:almostvalid}
d_{i}^{D}(j) \leq \delta \leq d_{i}^{D}(k).
\end{equation}

Almost valid distance drawings are not necessarily valid distance because there can be three vertices, $i \in V$, $j \in N^+(i)$ and $k \in N^-(i)$, such that 
$d_{i}^{D}(j) = \delta = d_{i}^{D}(k)$. Thus, condition \eqref{cond:valid} is not met for these vertices. It is also worth mentioning that a valid distance drawing 
may not be almost valid distance. But, we will show
that it is possible to obtain a valid distance drawing in $\mathcal{C}$ from an almost valid distance drawing. 

The $\argmax$ and $\argmin$ operators over a function give 
the elements of the domain of the function at which the function values are maximized and minimized, respectively. 
We use the $\argmax$ and $\argmin$ operators to define four important vertices for a vertex $a \in V$.
\begin{itemize}
\item
The \emph{farthest friend of $a$ on its left half of $\mathcal{C}$} is:
\[a_l^+:=\argmax_{\{j\colon\ j\in N^{+}[a],\ D(j)\in M_{l}(D(a))\}} 
d_{a}^{D}(j).\]
\item 
The \emph{closest enemy of $a$ on its left half of $\mathcal{C}$} is:
\[
a_l^-:=\argmin_{\{j\colon\ j\in N^{-}(a),\ D(j)\in M_{l}(D(a))\}} 
d_{a}^{D}(j).
\]
\item 
The \emph{farthest friend of $a$ on its right half of $\mathcal{C}$} is:
\[a_r^+:=\argmax_{\{j\colon\ j\in N^{+}[a],\ D(j)\in M_{r}(D(a))\}} 
d_{a}^{D}(j).\]
\item 
The \emph{closest enemy of $a$ on its right half of $\mathcal{C}$} is:
\[
a_r^-:=\argmin_{\{j\colon\ j\in N^{-}(a),\ D(j)\in M_{r}(D(a))\}} 
d_{a}^{D}(j).
\]
\end{itemize}
It is worth noting that $a_l^+$ and $a_r^+$ are well defined since, for all $a \in V$, $a\in N^+[a]$ and $D(a) \in M_l(D(a))\cap M_r(D(a))$.
On the other hand, $a_l^-$ and $a_r^-$ may not be defined since the set 
$\{j\colon\ j\in N^{-}(a)\land D(j)\in M_{r}(D(a))\}$ or the set 
$\{j\colon\ j\in N^{-}(a)\land D(j)\in M_{l}(D(a))\}$ may be empty. 
In that case, we say that the corresponding vertices, $a_l^-$ and/or $a_r^-$, do not exist. 

%Let $D$ be an almost valid distance drawing. We denote $\widehat{V}$ the set of vertices $a \in V$ such that $d(D(a),D(a_l^+)) = d(D(a),D(a_r^-))$ or $d(D(a),D(a_l^-)) = d(D(a),D(a_r^+))$. 

\begin{lemma}\label{lemma:almostvalidtovalid}
Let $S$ be a signed graph with an almost valid distance drawing in 
$\mathcal{C}$. Then, $S$ has a valid distance drawing in $\mathcal{C}$. 
\end{lemma}
\begin{proof}
Let $D$ be an almost valid distance drawing of $S$ in $\mathcal{C}$. 
We define $\widehat{V}\subseteq V$ as the set of all vertices for 
which restriction \eqref{cond:valid} is violated in $D$, i.e, 
\[\widehat{V}:=\{i \in V: \exists j\in N^{+}(i),\exists k \in N^{-}(i),\ 
d_{i}^{D}(j)=\delta=d_{i}^{D}(k)\}.\]
If $\widehat{V}=\emptyset$, then $D$ is valid distance. Hence, we assume 
that $\widehat{V}\neq \emptyset$.

For every $a \in \widehat{V}$,  $N^{-}(a)\neq \emptyset$. Moreover, since $D$ is an injection, 
for every  $a \in \widehat{V}$,
one, and only one, of the next situations occurs: 
\[d_{a}^{D}(a_r^+)=\delta=d_{a}^{D}(a_l^-), \mbox{ or } 
d_{a}^{D}(a_l^+)=\delta=d_{a}^{D}(a_r^-).\]
Otherwise, $D(a_l^+)=D(a_l^-)$ and $D(a_r^+)=D(a_r^-)$, which is not possible. 

Now, we consider a particular vertex $a \in \widehat{V}$.
Without loss of generality, we assume 
that $d_{a}^{D}(a_r^+)=\delta=d_{a}^{D}(a_l^-)$, (otherwise, by symmetry, we can reflect the drawing along the axis 
that goes through $a$ and the center of the circumference). 
Note that by the cyclic order  $a_l^+=(a_l^-+1)\bmod{n}$, and, if $a_r^-$ exists, $a_r^-=(a_r^++1)\bmod{n}$.

For every $b \in V \setminus \widehat{V}$, we define:
\[
\gamma_b := \min\{|d_{b}^{D}(b_l^-) - d_{b}^{D}(b_r^+)|, |d_{b}^{D}(b_r^-) - d_{b}^{D}(b_l^+)|\},
\]
where $d_{b}^{D}(b_l^-):=0$ if $b_l^-$ does not exist and $d_{b}^{D}(b_r^-):=0$ if $b_r^-$ does not exist. Since $b \in V \setminus \widehat{V}$, $\gamma_b >0$. We also define 
\[
\gamma := \left\{\begin{array}[l]{l l}
 \min_{b\in V\setminus \widehat{V}} \gamma_b    &   \mbox{ if }V \setminus \widehat{V} \neq \emptyset\\
2\pi     &\mbox{ if }V \setminus \widehat{V} = \emptyset. 
\end{array}\right.
\]

We transform $D$ into an almost valid distance drawing $D'$ for which the set of  
vertices that violate restriction \eqref{cond:valid} does not contain $a$
and does not contain any new vertex.
We use $\bar{a}:={(D(a)+\pi)}\bmod{2\pi}$ to denote the antipodal point of $D(a)$. Now, we define the following value: 
\[\epsilon:=\frac{1}{4}\min\{\min_{i\in V} 
d_{i}^{D}\left({(i+1)}\bmod{n}\right), d_{\bar{a}}^{D}\left(a_l^-\right), \gamma\}.\]
Since $D$ is an injection, $\min_{i\in V} 
d_{i}^{D}\left((i+1)\bmod{n}\right) > 0$. Since $a \in \widehat{V}$, $a_l^-$ cannot be the antipodal point of $a$, thus,  $d_{\bar{a}}^{D}\left(a_l^-\right) > 0$.
 Also, by definition $\gamma> 0$, thus, $\epsilon > 0$. 
Furthermore, due to the definition of $\epsilon$ it also holds that 
$\delta > \epsilon$.
Now, we define the following drawing $D':V\rightarrow \mathcal{C}$,
\[
D'(i)=\left\{ 
\begin{array}
[l]{l l}
D(i) &\mbox{ if }i\neq a,\\
D(i)+\epsilon  &\mbox{ if }i= a.
\end{array}
\right. 
\]

We first observe that, since $\epsilon<d_{i}^{D}\left((i+1)\bmod{n}\right)$, the 
cyclic order of $V$ induced by $D'$ is equal to the cyclic order of $V$ induced by $D$. 
Therefore, the cyclic labeling according to $D$ is the same as the cyclic
labeling according to $D'$. 
Moreover, since $\epsilon< d_{\bar{a}}^{D}\left(a_l^-\right)$, in this new drawing the vertices $a_l^+$, $a_r^+$, $a_l^-$ and 
$a_r^-$ remain the same.

We show now that $a$ does not violate restriction  \eqref{cond:valid} 
in $D'$. We analyze the distance from $a$ to its farthest friends
and closest enemies in $D'$. We start with the distances between $a$ and its farthest friends in $D'$. 
\begin{eqnarray*}
d_{a}^{D'}(a_r^+)&=&d_{a}^{D}(a_r^+)-\epsilon=\delta-\epsilon <\delta,
\end{eqnarray*}
and
\begin{eqnarray*}
d_{a}^{D'}(a_l^+)&=&d_{a}^{D}(a_l^+)+\epsilon\\
&=&d_{a}^{D}(a_l^-)-d_{a_l^-}^{D}(a_l^+)+\epsilon\\
&=&\delta-d_{a_l^-}^{D}(a_l^+)+\epsilon\\ &<&\delta,
\end{eqnarray*}
where the last inequality is obtained since $a_l^+=(a_l^-+1)\bmod{n}$, and hence
\[d_{a_l^-}^{D}(a_l^+)=d_{a_l^-}^{D}\left((a_l^-+1)\bmod{n}\right)>\epsilon.\]

On the other hand, if we repeat the analysis for the distances between $a$ and its closest enemies in $D'$, we obtain:
\begin{eqnarray*}
d_{a}^{D'}(a_l^-)&=& d_{a}^{D}(a_l^-)+\epsilon
=\delta+\epsilon
> \delta,
\end{eqnarray*}
and
\begin{eqnarray*}
d_{a}^{D'}(a_r^-)&=& d_{a}^{D}(a_r^-)-\epsilon\\
&=&d_{a}^{D}(a_r^+)+d_{a_r^+}^{D}(a_r^-)-\epsilon\\
&=& \delta+d_{a_r^+}^{D}(a_r^-)-\epsilon\\
&>& \delta,
\end{eqnarray*}
where the last inequality is obtained since $a_r^-=(a_r^++1)\bmod{n}$, and hence 
\[d_{a_r^+}^{D}(a_r^-)= d_{a_r^+}^{D}\left((a_r^++1)\bmod{n}\right)> \epsilon.\]
Therefore, vertex $a$ does not violate restriction \eqref{cond:valid}
in $D'$. 

Consider a vertex $b\neq a$. 
The distance between $b$ and any other vertex $c\neq a$ is the same in $D$ and in $D'$. Therefore, if none of the vertices $b_l^+$, $b_r^+$, $b_l^-$ or $b_r^-$ is equal to $a$,  $b$ belongs to $\widehat{V}$ or meets condition \eqref{cond:valid} in $D'$ as it did in $D$. 
%\begin{itemize}

If $a=b_l^+$, then $b \in N^+(a)$ and $D(b) \in M_r(D(a))$, thus: % $d_{b}^{D'}(b_l^+)=d_{b}^{D'}(a)$ and:
\[
d_{b}^{D'}(b_l^+)=d_{b}^{D'}(a)= d_{b}^{D}(a) - \epsilon < d_{b}^{D}(a) \leq d_{b}^{D}(b_r^-)
= d_{b}^{D'}(b_r^-).
\]

If $a= b_r^-$, then  $b \in N^-(a)$ and $D(b) \in M_l(D(a))$, thus: % $d_{b}^{D'}(b_r^-)=d_{b}^{D'}(a)$ and:
\[
d_{b}^{D'}(b_r^-)=d_{b}^{D'}(a) = d_{b}^{D}(a) + \epsilon > d_{b}^{D}(a) \geq d_{b}^{D}(b_l^+)= d_{b}^{D'}(b_l^+).
\]

If now $b \in V\setminus \widehat{V}$ and $a=b_r^+$, then $b \in N^+(a)$ and $D(b) \in M_l(D(a))$, thus: %$d_{b}^{D'}(b_r^+)=d_{b}^{D'}(a)$ and:
\[
d_{b}^{D'}(b_r^+)=d_{b}^{D'}(a) = d_{b}^{D}(a) + \epsilon < d_{b}^{D}(a) + \gamma_b \leq d_{b}^{D}(b_l^-) = d_{b}^{D'}(b_l^-).
\]

And if $b \in V\setminus \widehat{V}$ and $a= b_l^-$, then  $b \in N^-(a)$ and $D(b) \in M_r(D(a))$, thus:% $d_{b}^{D'}(b_l^-)=d_{b}^{D'}(a)$ and:
\[
d_{b}^{D'}(b_l^-)=d_{b}^{D'}(a) = d_{b}^{D}(a) - \epsilon > d_{b}^{D}(a) - \gamma_b \geq  d_{b}^{D}(b_r^+)= d_{b}^{D'}(b_r^+).
\]

To conclude the proof, we  need to see that there is no vertex $b \in \widehat{V}$ such that $a= b_r^+$ or $a= b_l^-$. If there is a vertex $b \in \widehat{V}$ such that $a= b_r^+$, then $D(b) \in M_l(D(a))$, $b \in N^+(a)$, and $d_{a}^{D}(b)=\delta$. Since $a_l^- \in N^-(a)$, $b\neq a_l^-$.
Since $D(a_l^-)\in M_l(D(a))$ and $d_{a}^{D}(a_l^-)=\delta$, $D(b)=D(a_l^-)$. Which contradicts the fact that $D$ is an injection. A symmetric argument allows us to conclude that there is no vertex $b \in \widehat{V}$ such that $a= b_l^-$. 

In conclusion, converting $D$ into $D'$ has decreased the size of $\widehat{V}$ by at least one. Hence, if we repeat this process at most $|\widehat{V}|$ times, we obtain a valid distance drawing for $S$ in $\mathcal{C}$.
\end{proof}

\subsection{Circular-Arc Graphs}
A \emph{circular-arc} graph is the intersection graph of a set of arcs on the circumference. It has one vertex for each arc in the set and an edge between every pair of vertices corresponding to arcs that intersect. The set of arcs that corresponds to a graph $G$ is called a \emph{circular-arc model} of $G$. A \emph{proper circular-arc} graph is a circular-arc graph for which there exists a corresponding circular-arc model such that no arc properly contains another. Such model is called  \emph{proper circular-arc} model. On the other hand, a \emph{unit circular-arc} graph is a circular-arc graph for which there exists a corresponding circular-arc model such that each arc is of equal length. Such model is called \emph{unit circular-arc} model.    

The first characterization of circular-arc graphs is due to Alan Tucker in \cite{tucker1970characterizing}. The same author in \cite{doi:10.1137/0209001} also presented the first polynomial-time recognition algorithm for this family of graphs, which runs in $\mathcal{O}(n^3)$ time. Ross M. McConnell in \cite{mcconnell2003linear} presented the first $\mathcal{O}(n+m)$ time recognition algorithm for circular-arc graphs. Recognition of a proper circular-arc graph and construction of a proper circular-arc model can both be performed in time $\mathcal{O}(n+m)$ as it was proved by Deng et al. in \cite{deng1996linear}. On the other hand, Duran et al. in \cite{duran2006polynomial} presented a $\mathcal{O}(n^2)$ time algorithm for the recognition of unit circular-arc graphs.  Later, Lin and Szwarcfiter gave an $\mathcal{O}(n)$ time recognition algorithm for unit circular-arc graphs that also constructs a unit circular-arc model for the graph in consideration. 

From the definition, we can see that every proper circular-arc graph is also a circular-arc graph. Nevertheless, the opposite contention does not hold. For example the \emph{net} (see Figure \ref{fig}) is a circular-arc graph that does not have a proper circular-arc model. In the same line, every unit circular-arc graph is also a proper circular-arc graph. However, Alan Tucker in \cite{tucker1974structure} gave a characterization of proper circular-arc graphs, which are not unit circular-arc graphs. It is worth noting that Tucker's characterization uses crucially the fact that all the unit arcs are closed, or all are open. Kaplan and Nassbaum in \cite{kaplan2009certifying} pointed out that the family of unit circular-arc graphs does not change if all the arcs are restricted to be open or all are restricted to be closed.
Nevertheless, the family changes if we allow unit circular-arc models to have arcs that are open and closed. Every proper circular-arc graph has an arc model with arcs of the same length, where arcs can be open and closed. This fact was pointed out in \cite{kaplan2009certifying} and deduced from a construction presented in \cite{tucker1974structure}. Although, none of these two articles stated this fact as a result. Since we use it later, we believe it is important to express it as a lemma. 

For the sake of completeness, we also point out the fact that circular-arc models and proper circular-arc models are not restricted to have closed or open arcs. Indeed, it is always possible to perturb the arcs in a circular-arc model and a proper circular-arc model so that no two arcs share an endpoint. Hence, whether arcs are open or closed in these models does not affect the corresponding graph. 
 
\begin{lemma}{(\cite{tucker1974structure})}\label{lemma:unitmodelproper}
Let $G$ be a proper circular-arc graph. Then $G$ has an arc model 
with 
arcs of the same length, where arcs can be open and closed.  
 \end{lemma}

The construction of an arc model with arcs of the same length, where arcs can be open and closed, for a proper circular-arc graph is the same construction described in the proof of Theorem 4.3 in \cite{tucker1974structure}. However, when unit arcs create new intersections in their endpoints, and they cannot be moved to avoid these new intersections, the arc's endpoint is deleted from the arc, creating an open arc of unit length.

\section{Signed Graphs with a Valid Distance Drawing in the 
Circumference}\label{sec:characterization}

Given a signed graph $S$, a \emph{completion} of $S$ is a set of decisions of the type $ij \in E^+$ or $ij \in E^-$ for all pairs $ij \notin E^+ \cup E^-$. Given a signed graph $S$, we say that $C_S$, a complete signed graph on the same set of vertices than $S$, is a \emph{completion} of $S$ if $S$ is a signed subgraph of $C_S$. We use $C_S^+$ to denote the positive subgraph of $C_S$.

\begin{lemma}\label{thm:charact}
Let $S$ be a signed graph, and $C_S$ be a completion of $S$ such that $C_S^+$ is a proper
circular-arc graph. Then $S$ has a valid distance drawing in $\mathcal{C}$  
\end{lemma}
\begin{proof}
%Let us first point out the fact that, if a signed graph $S$ has a 
%valid distance
%drawing in $\mathcal{C}$, then any signed subgraph $H$ of $S$ also 
%has a 
%valid distance drawing in $\mathcal{C}$. 
%This affirmation can be obtained by simply considering exactly the
%same 
%valid distance drawing $D$ for $S$ but now for $H$. Indeed, all the
%restrictions 
%(\ref{cond:valid}) that need to be satisfied by a 
%valid distance drawing for $H$ are already satisfied by $D$.
%Moreover, by deleting
%edges or vertices from $S$ to obtain $H$, the only thing that may 
%happen 
%is that some of the restrictions (\ref{cond:valid})
%that are satisfied by $D$ are not required by $H$, but this does not
%harm
%the validity of $D$ for $H$. 
 
Consider a signed graph $S$ with a completion $C_S$ such that $C_S^+$ 
is a proper circular-arc graph. 
%We construct a valid distance drawing for $S$ in $\mathcal{C}$. 
%Hence, due to
%the affirmation of the previous paragraph, this valid distance 
%drawing
%is also a valid distance drawing for $G$ (since $G$ is a signed 
%subgraph of $C_G^+$).
%
Lemma \ref{lemma:unitmodelproper} says that $C_S^+$ has an arc model
with
arcs of the same length, let say of length $\delta$. Let $s_i$ be 
the anticlockwise end of the arc corresponding to vertex $i$
in a such model. We define the drawing 
$D:V\rightarrow \mathcal{C}$ as
$D(i):=s_i$. The drawing $D$ satisfies $d_{i}^{D}(j)\leq \delta$ for
all 
$ij \in E^+$ and $d_{i}^{D}(j)\geq \delta$ for all $ij \in E^-$.
Hence, 
$D$ is almost valid distance. Therefore, by Lemma 
\ref{lemma:almostvalidtovalid}, $D$ can be transformed into a valid 
distance drawing of $S$ in $\mathcal{C}$.
\end{proof}

\begin{figure}[t]
 \begin{minipage}[t]{.49\linewidth}
 \centering\includegraphics[scale=0.8]{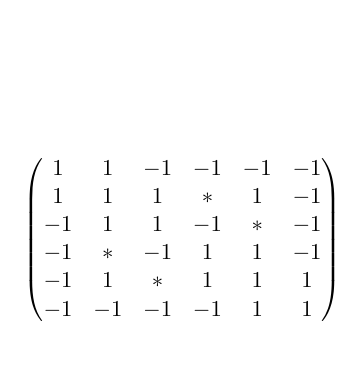}
 \subcaption{Augmented adjacency matrix of the signed graph $S$.}\label{fig:3:4}
 \end{minipage}
 \begin{minipage}[t]{.49\linewidth}
 \centering
\includegraphics[scale=0.8]{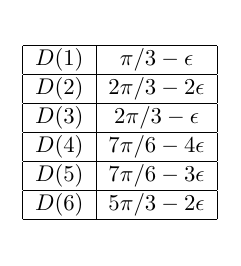}
\subcaption{Valid distance drawing for $S$ in $C$.}\label{fig:3:5}
 \end{minipage}
\begin{minipage}[t]{.32\linewidth}
 \centering\includegraphics[scale=0.6]{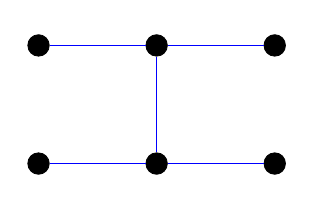}
 \subcaption{Positive subgraph of $S$, and positive subgraph of the completion that gives $-1$ to the two missing edges.}\label{fig:3:1}
 \end{minipage}
 \begin{minipage}[t]{.32\linewidth}
 \centering\includegraphics[scale=0.6]{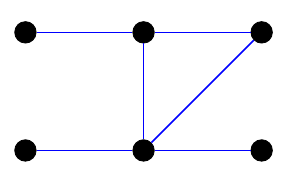}
 \subcaption{Positive subgraph of the two completions of $S$ that give $1$ to one edge and $-1$ to the other. }\label{fig:3:2}
 \end{minipage}
 \begin{minipage}[t]{.32\linewidth}
 \centering\includegraphics[scale=0.6]{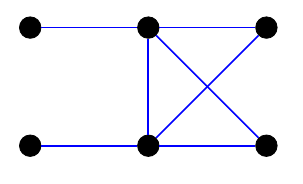}
 \subcaption{Positive subgraph of the completion of $S$ that gives $1$ to the two missing edges.}\label{fig:3:3}
 \end{minipage}

 \caption{This figure shows a signed graph $S$ via its augmented adjacency matrix in Subfigure \ref{fig:3:4}, with a valid distance drawing that is shown in Subfigure \ref{fig:3:5}. Nevertheless, 
 $S$ does not have a completion whose positive subgraph is proper circular-arc. Subfigures \ref{fig:3:1}, \ref{fig:3:2}, and \ref{fig:3:3} show the positive subgraphs that are obtained under  any completion of $S$. None of these graphs is proper circular-arc since all of them contain as induced subgraph the \emph{claw}, defined as the complement of $C_3 \cup K_1$. 
 \label{fig:3}}
 \end{figure}

The opposite direction of Lemma \ref{thm:charact} is not true. Figure \ref{fig:3} shows a signed graph with a valid distance drawing in $\mathcal{C}$ and for which no completion produces a proper circular-arc graph as positive subgraph. We show a weaker result instead. We show that if a signed graph has a valid distance drawing in the circumference, it has a completion such that its positive subgraph is circular-arc. We use the matrix characterization for circular-arc graphs given by Alan Tucker in \cite{tucker1971matrix}.

For any graph $G$, the \emph{adjacency matrix} $A(G)$ of $G$ is the symmetric $n\times n$ matrix with entry $(i,j)$ equal to $1$ if $ij \in E$ and equal to $0$ if $ij \notin E$. The \emph{augmented adjacency matrix} of $G$ is the matrix $A^*(G):=A(G)+I$, where $I$ denotes the $n \times n$ identity matrix. On the other hand, let $A$ be a symmetric $(0,1)$-matrix with $1$'s on the main diagonal. Let $D_i$ be the set of $1$'s in column $i$ starting at the main diagonal and going down and around until the first $0$ is encountered. Let $R_i$ be the set of $1$'s in row $i$ starting at the main diagonal and going right and around until the first $0$ is encountered. Then $A$ has the \emph{quasi-circular} $1$'s property if the union of all $D_i$'s and $R_i$'s contain all $1$'s in $A$. Alan Tucker in \cite{tucker1971matrix} proved the following theorem. 
\begin{theorem}[Theorem 2 in \cite{tucker1971matrix}]\label{thm:tuckercharac}
A graph $G$ is circular-arc if and only if its vertices can be indexed so that $A^*(G)$ has the quasi-circular $1$'s property. 
\end{theorem}

We consider a natural extensions for the definitions of adjacency matrix, augmented adjacency matrix, and quasi-circular $1$'s property in the context of signed graphs. Let $S$ be a signed graph, the \emph{adjacency matrix} of $S$ is the symmetric $n \times n$ matrix $A(S)$ with entry $(i,j)$ equal to $1$ if $ij \in E^+$, equal to $-1$ if $ij \in E^-$, and equal to $*$ if $ij \notin E^+ \cup E^-$. The \emph{augmented adjacency matrix} of a signed graph $S$ is the matrix $A^*(S) := A(S) + I$. On the other hand, let $A$ be a symmetric matrix with entries taking values in $\{-1,*,1\}$ and with $1$'s in the main diagonal. Let $D_i$ be the set of $1$'s in column $i$ starting at the main diagonal and going down and around until the first $-1$ is encountered. Similarly, let $R_i$ be the set of $1$'s in row $i$ starting in the main diagonal and going right and around until the first $-1$ is encountered. Matrix $A$ is said to have the \emph{quasi-circular} $1$'s property if the union of all $D_i$'s and $R_i$'s contains all the $1$'s in $A$.

In the context of matrices with entries in $\{-1,*,1\}$, a completion is a set of decisions that transform each $*$ into either a $1$ or a $-1$.
\begin{lemma}\label{lem:quasicirc-compquasicirc}
Let $A$ be a symmetric matrix with entries taking value in the set $\{-1,*,1\}$ and with $1$'s in the main diagonal.  The matrix $A$ has the quasi-circular $1$'s property if and only if it has a completion with the quasi-circular $1$'s property. 
\end{lemma}
\begin{proof}
 On one hand, if a matrix $A$ has the quasi-circular $1$'s property, then, transforming all $*$'s in the union of the $D_i$'s and the $R_i$'s into $1$'s, and all $*$'s outside that union into $-1$, we obtain a completion of $A$ with the quasi-circular $1$'s property.
 
On the other hand, if a completion of the matrix $A$ has the quasi-circular $1$'s property, $A$ inherits the property.
To see this, 
%we transform that completion back into $A$. 
%In that transformation, some $-1$'s and $1$'s are transformed into $*$,
%but a $-1$ or a $1$ is never transformed into a $1$ or a $-1$, 
%respectively. 
%Therefore, if we 
assume that $A$ does not have the quasi-circular $1$'s property but it has a completion with the quasi-circular $1$'s property. Hence, $A$ has an entry $(j,l)$ equal to $1$, an entry $(j,k)$ equal to $-1$, and an entry $(i,l)$ equal to $-1$ such that $i<j<k<l$.  These three entries maintain their value in any completion of $A$. Consequently, no completion can have the quasi-circular $1$'s property, which is a contradiction.
\end{proof}

In terms of signed graphs, Lemma \ref{lem:quasicirc-compquasicirc} and Theorem \ref{thm:tuckercharac} allow us to conclude the following corollary. 
\begin{corollary}\label{cor:signedqciffcompletsignedqc}
Let $S$ be a signed graph. The matrix $A^*(S)$ has the quasi-circular $1$'s property if and only if $S$ has a completion $C_S$ such that $C_S^+$ is circular-arc.
\end{corollary}

Now, using Lemma \ref{lem:quasicirc-compquasicirc} and Corollary \ref{cor:signedqciffcompletsignedqc}, we show that if a signed
graph has a valid distance drawing in the circumference, it has a completion
such that its positive subgraph is circular-arc. 

\begin{lemma}\label{lem:vddtoqc1matrix}
Let $S$ be a signed graph with a valid distance drawing in $\mathcal{C}$. 
Then, $S$ has a completion $C_S$ such that its positive subgraph $C_S^+$ is 
circular-arc.
\end{lemma}
\begin{proof}
Let $S=(V,E^+\cup E^-)$ be a signed graph with a valid distance drawing in $\mathcal{C}$. Let $D$ be a such drawing. We show that the ordering of $V$ induced by $D$ indexes $V$ such that $A^*(S)$ has the quasi-circular $1$'s property. By contradiction, assume that this is not the case. Hence, $A^*(S)$ has an
 entry $(j,l)$ equal to $1$ that is not contained neither in $D_l$ nor in $R_j$. Without loss of generality, assume that $j<l$. Since neither $D_l$ nor $R_j$ contain the entry $(j,l)$, there is a $-1$ in a position $(j,k)$ of $A^*(S)$ such that $j<k<l$, and a $-1$ in a position $(i,l)$ of $A^*(S)$ such that $i<j<l$. Therefore, $\{il,jk\} \subseteq E^-$ and $jl \in E^+$.

There are two possibilities for $D(l)$, either $D(l)$ is in the right half of $D(j)$ in $\mathcal{C}$
or $D(l)$ is in the left half of $D(j)$. 
If $D(l)$ is in the left half of $D(j)$, $D(k)$ is also in the left half of $D(j)$ in $\mathcal{C}$, since $j<k<l$ according to $D$. Moreover, 
\[
d_{j}^{D}(k) < d_{j}^{D}(l),
\]
which contradicts the fact that $D$ is valid distance because $jk \in E^-$ and $jl \in E^+$. Therefore, $D(l)$ is in the right half of $D(j)$ in $\mathcal{C}$.
Since $i<j<l$ according to $D$ (which in circular terms means $j<l<1<i<j$), $D(i)$ is also in the right half of $D(j)$. Furthermore, 
\[
d_{j}^{D}(i) < d_{j}^{D},
\]
which also contradicts the fact that $D$ is valid distance since $ji \in E^-$ and $jl \in E^+$.

Hence, such triplet of entries that brake the quasi-circular $1$'s  property cannot exist. Therefore, $A^*(S)$ has the quasi-circular $1$'s property, and  Lemma \ref{lem:quasicirc-compquasicirc}, and Corollary \ref{cor:signedqciffcompletsignedqc} allow us to conclude that $S$ has a completion whose positive subgraph is circular-arc.  
\end{proof}

Let $A$ be a symmetric matrix with entries in $\{-1,*,1\}$ and with $1$'s on the main diagonal.  
Let $L_i$ be the set of $1$'s in row $i$ starting at the main diagonal and going left and around until the first $-1$ is encountered. Let $R_i$ be the set of $1$'s in row $i$ starting at the main diagonal and going right and around until the first $-1$ is encountered. Then $A$ has the \emph{circular} $1$'s property if the union of all $L_i$'s and $R_i$'s contain all $1$'s in $A$. 
That is, $1$'s in each row of $A$ (and columns due to symmetry) appear in a circular fashion, potentially broken by $*$, but never broken by $-1$.
It is worth mentioning that a symmetric matrix $A$ with entries in $\{-1,*,1\}$ and with $1$'s on the main diagonal that has the circular $1$'s property, also has a completion that has the circular $1$'s property. It is enough to transform all $*$ in the union of $L_i$'s and $R_i$'s into $1$'s and 
those outside that union into $-1$'s. Nevertheless, this completion 
may not be symmetric. Therefore, it does not translate directly to completions of signed graphs. 

We say that a symmetric matrix $A$ with entries in $\{-1,*,1\}$ and with $1$'s on the main diagonal has the \emph{circularly compatible} $1$'s property if $A$ has the circular $1$'s property, and if, after inverting and/or cyclically permuting the order of rows and corresponding columns of $A$ the following holds: 
let $j$ be the smallest index such that entry $(j,1)$ is equal to $-1$, 
therefore, for all $i\in\{2,\ldots,j-1\}$ if entry $(i,2)$ is equal to $-1$, then entry $(i,1)$ equals $*$, unless one of these columns is all $1$'s or all $-1$'s (except in the main diagonal).
 
It is noteworthy that the definition of the circularly compatible $1$'s property for matrices with entries in the set $\{-1,*,1\}$, but restricted to the case when no entry is equal to $*$, says that a such matrix has the  circular $1$'s property and also, after inverting and/or cyclically permuting the order of rows and corresponding columns, the last $1$ of the circular set in the second column is always at least as low as the last $1$ of the circular set in the first column, unless one of these columns is all $1$'s or all $-1$'s.
On the other hand, not every matrix with entries in the set $\{-1,*,1\}$ that has the circularly compatible $1$'s property has a completion with the circularly compatible $1$'s property. The matrix in Subfigure \ref{fig:3:4} is an example of this.

\begin{lemma}\label{lem:vddacc1p}
Let $S$ be a signed graph, and $D$ be a valid distance drawing of $S$ in $\mathcal{C}$. Then, $A^*(S)$ has the circularly compatible $1$'s property when the vertices of $S$ are indexed according to the order induced by $D$.
\end{lemma}
\begin{proof}
Let $S$ be a signed graph and $D$ be a valid distance drawing of $S$ in $\mathcal{C}$. Let $A^*(S)$ be the augmented adjacency matrix of $S$ with the vertices of $S$ labeled according to the order induced by $D$.

By contradiction, assume  that $A^*(S)$ hasn't the circularly compatible $1$'s property because it does not have the circular $1$'s property. Therefore, without loss of generality (by cyclically permuting the order of rows and corresponding columns), assume that there are columns $i<j<k$ such that entry $(1,i)$ is equal to $-1$, entry $(1,j)$ is equal to $1$ and entry $(1,k)$ is equal to $-1$ in $A^*(S)$. Hence, $\{1i,1k\} \subseteq E^-$ and $1j \in E^+$. Since  $i<j<k$ according to $D$, then $d_{1}^{D}(i)< d_{1}^{D}(j)$ or $d_{1}^{D}(k)< d_{1}^{D}(j)$. In any of these two cases, the necessary condition to be a valid distance drawing is broken. Which is a contradiction since $D$ is a valid distance drawing. 

Now assume  that $A^*(S)$ hasn't the circularly compatible $1$'s property because there are entries $(i,2)$ equal to $-1$, $(i,1)$ equal to $1$ and $(j,1)$ equal to $-1$ such that $2<i<j$. Therefore, $\{i2,j1\} \subseteq E^-$ and $i1 \in E^+$. There are two options for $i$, $D(i)$ is in the left half of $D(1)$ in $\mathcal{C}$, or $D(i)$ is in the right half of $D(1)$ in $\mathcal{C}$.  

If $D(i)$ is in the left half of $D(1)$ in $\mathcal{C}$, since $1<2<i$, $D(2)$ is also in the left half of $D(1)$ in $\mathcal{C}$ and 
$d_{i}^{D}(2) < d_{i}^{D}(1)$.  Which is a contradiction with the fact that $D$ is valid distance. 
Therefore, $D(i)$ is in the right half of $D(1)$ in $\mathcal{C}$. Since $i<j$, then $D(j)$ is also in the right half of $D(1)$ in $\mathcal{C}$. Moreover, 
$d_{1}^{D}(j)< d_{1}^{D}(i)$, which is a contradiction with the fact that $D$
is valid distance. 
In conclusion, $A^*(S)$ has the  circularly compatible $1$'s property.
\end{proof}

Lemma \ref{thm:charact} and Lemma \ref{lem:vddacc1p} allow us to state the following theorem. 
\begin{theorem}\label{thm:characcomplete}
A complete signed graph $S$ has a valid distance drawing in $\mathcal{C}$ if and only if $S^+$ is a proper circular-arc graph.  
\end{theorem}
\begin{proof}
 Let $S$ be a complete signed graph. On the one hand, Lemma \ref{thm:charact} says that if $S^+$ is a proper circular-arc graph, then $S$ has a valid distance drawing in $\mathcal{C}$. On the other hand, Lemma \ref{lem:vddacc1p} says that if $S$ has a valid distance drawing in $\mathcal{C}$, then $A^*(S)$ has the circularly compatible $1$'s property. This lemma, together with Theorem 6 in \cite{tucker1971matrix}, that says that a graph is proper circular-arc if and only if its augmented adjacency matrix has the circularly compatible $1$'s property, complete the proof. 
\end{proof}

\section{Concluding Remarks}\label{sec:concludingremarks}
We would like to finish this document by mentioning two exciting problems left open in this document. First of all, a characterization of the family of signed graphs (not necessarily complete) with a valid distance drawing in $\mathcal{C}$ is still missing. In Lemma \ref{lem:vddacc1p}, we give a necessary condition for a graph to belong to that family. Is that condition sufficient?. Finally, and undoubtedly related to the previous question,  it is still an open problem to find out whether the SCFE problem in the circumference is NP-Complete or not. This problem is NP-Hard since there is a reduction from the SCFE on $\mathbb R$. Let $S$ be an instance of SCFE on $\mathbb R$ and let  $S'$ be the signed graph $S$ with an additional vertex $x$  that is an enemy of all other vertices. Then, a valid distance drawing of $S$ in  $\mathbb R$ induces a  valid distance drawing of $S'$ in  $\mathcal C$ (insert the drawing of $S$ on a quarter of the circumference and place $x$ in the opposite side). Conversely,  a valid distance drawing of $S'$ in  $\mathcal C$  induces a valid distance drawing of $S$ in  $\mathbb R$  (cut the circumference at the position of $x$ and straighten it). Finding a characterization of the family of signed graphs with a valid distance drawing in $\mathcal{C}$ certainly would help answer this last question.

%\bibliographystyle{plainurl}% the mandatory bibstyle
%\bibliography{references}

\end{document}